\documentclass[11pt]{article}
\usepackage[utf8]{inputenc}
\usepackage[margin=2.3cm]{geometry}
\usepackage{amsmath,amsthm,amsfonts,amssymb}
\usepackage{graphicx,latexsym}
\usepackage{graphicx}
\usepackage[figuresright]{rotating}  
\usepackage{lscape}
\usepackage{cite}
\usepackage{float}

\usepackage{url} 
\usepackage[ruled,lined]{algorithm2e}
\usepackage{graphicx}
\usepackage{tikz}
\newtheorem{theorem}{Theorem}
\newtheorem{definition}{Definition}
\usetikzlibrary{automata, positioning, arrows}
\usepackage{cite}
\title{Pancyclicity in the  Cartesian Product $(K_9-C_9 )^n$}
\author{Syeda Afiya\\
School of Advanced Sciences\\
Vellore Institute of Technology, Chennai-600127\\ India \\ 
afiyazaakia@gmail.com
\and M. Rajesh\\
School of Computer Science and Engineering \\ 
Vellore Institute of Technology, Chennai-600127\\ India\\
 rajesh.m@vit.ac.in}
\date{June 2022}
\begin{document}
\maketitle
\begin{abstract}
  A graph $G$ on $m$ vertices is pancyclic if it contains cycles of length $l$, $3\leq l \leq m$ as subgraphs in $G$. The complete graph $K_{9}$ on 9 vertices with a cycle $C_{9}$ of length 9 deleted from $K_{9}$ is denoted by $(K_{9}-C_{9})$. In this paper, we prove that $(K_{9}-C_{9})^{n}$, the Cartesian product of $(K_{9}-C_{9})$ taken $n$ times, is pancyclic. 
\end{abstract}

{\bf{Keywords}} pancyclicity, cartesian product, $(K_9-C_9 )^n$,  hamiltonian

\section{Introduction}
The concept of pancyclicity was first proposed by Bondy {\cite{bondy1971pancyclic}}. Pancyclicity in a network
enables the study of algorithms designed for cycles in a network. A cycle $C_l$ of length $l$ is called an $l$-cycle. A graph $G$ is said to be pancyclic if it has an $l$-cycle for every $l$ between $3$ and $n$ where $n$ is the number vertices in the graph $G$. Analogously, $G$ is edge-pancyclic if for every edge $e$ of $G$ and every integer $n$, $3\leq n \leq |V(G)|$, there is a  cycle of length  in $G$ using $e$ \cite{araki2003edge}. The circumference $c(G)$ of a graph $G$ is the length of a longest cycle contained in $G$. If $c(G)=n$, then the graph is hamiltonian. Pancyclicity has important applications in  parallel processing systems, storage schemes of logical data structures, layout of circuits in VLSI and designing low cost simple algorithms \cite{fan2005node,araki2002pancyclicity}.\\
\indent Several parallel algorithms depends on circle structures \cite{leighton2014introduction} and pancyclicity determines whether the topology of the system admits cycles of all lengths. To execute algorithms corresponding to cycle structure in an interconnection network, any two processes that are adjacent in the cycle structure are mapped to two adjacent nodes of the network \cite{li2013vertex}. Hence, a well-organized mapping requires that the system owns a cycle of a specified length. Further, these cycles may be used as control structures for distributed systems. Thus it is desirable for interconnection networks to be pancyclic. Crossed cubes, Möbius cube, $k$-Ary $n$-Cube, prism, WK-Recursive network, circulant graphs and OTIS-mesh network have been proved to be pancyclic \cite{fang2007m,lin2011panconnectivity,malekimajd2011pancyclicity,goddard2001pancyclicity, shafiei2011pancyclicity,bogdanowicz1996pancyclicity,yang2006pancyclicity}.
\section{Preliminaries}
In 1971 Bondy \cite{bondy1971pancyclic} posed the following meta conjecture:
Almost any nontrivial condition on a graph which implies that the graph is hamiltonian also implies that it is pancyclic (except for maybe a simple family of exceptional graphs).
For several sufficient conditions Bondy's meta conjecture has been verified. This is the motivation to examine these sufficient conditions even for vertex pancyclicity.  The first sufficient condition for a graph to be hamiltonian is due to Dirac in 1952.\\

\begin{theorem} {\bf \cite{dirac1952some}}
If $G$ is a graph of order $n\geq 3$ such that $\delta(G)\geq n / 2$ , then $G$ is hamiltonian.
\end{theorem}
\begin{theorem} {\bf \cite{ore1960note}}
 Let $G$ be a graph with $n$ vertices and let $u$ and $v$ be non adjacent vertices in $G$, such that $d(u)+d(v) \geq n$. Then $G$ is hamiltonian.
\end{theorem}
\indent In this paper, we study the pancyclicity of the graph $(K_9-C_9)^n$, $n\geq 1$. We have the following result on pancyclicity.
\begin{theorem} {\bf \cite{bondy1971pancyclic}}
Let $G$ be hamiltonian and suppose that $|E(G)| \geq \frac{n^2}{4}$, where $n=|V(G)|$. Then $G$ is either pancyclic or else is the complete bipartite graph $K_{\frac{n}{2},\frac{n}{2}}$.
\end{theorem}
\indent However, a graph $G$ on $n$ vertices is pancyclic does not necessarily imply that $|E(G)| \geq \frac{n^2}{4}$. For example, the circulant graph $C(10,1,2)$ on 10 vertices and 20 edges is pancyclic. See Figure 1 (a) and (b).
\begin{figure}[h]
    \centering
    \includegraphics[width=13cm]{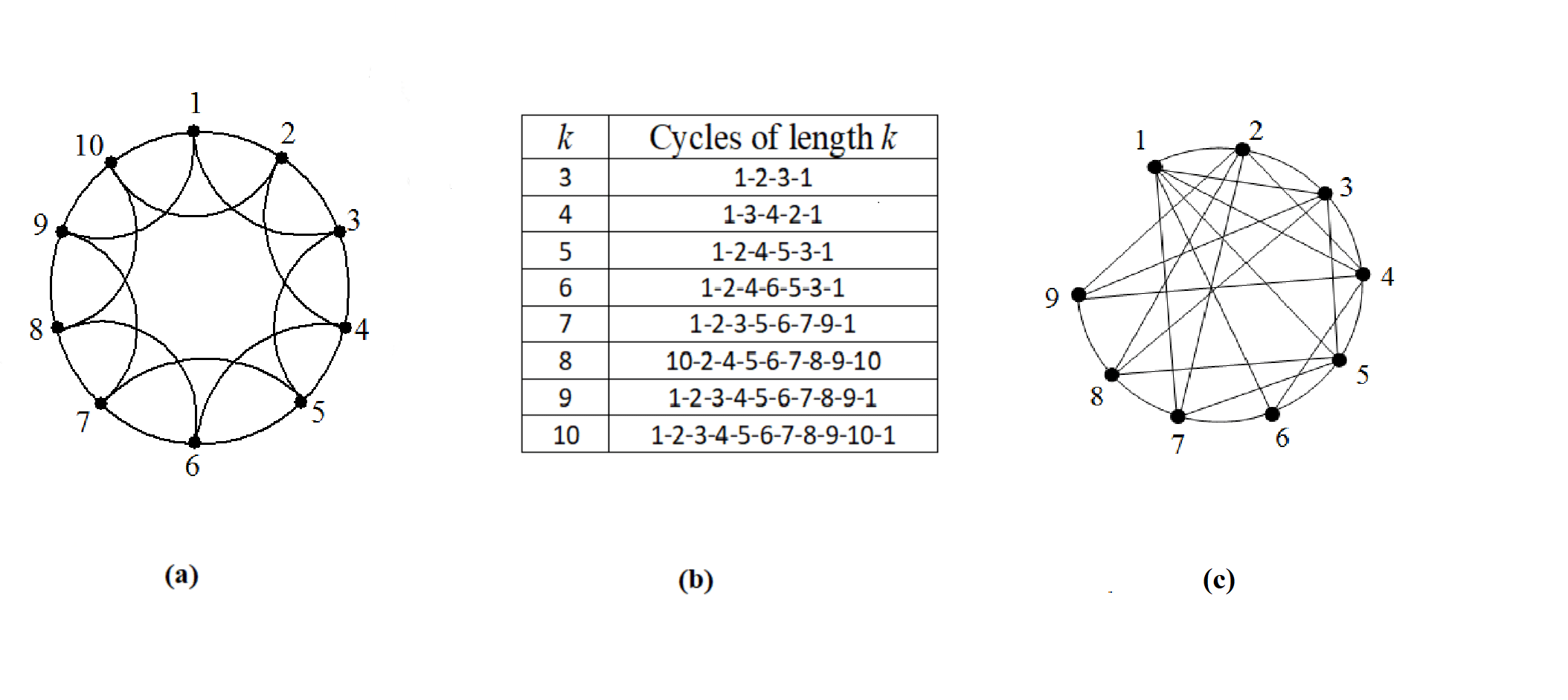}
    \caption{(a) Circulant $C(10,1,2)$ (b) Cycles of length $k$ (c) $(K_9-C_9)$ with a hamiltonian cycle $1-3-4-5-6-7-8-9-2-1$}
    \label{Figure 1}
\end{figure}\\
\begin{definition}{\bf \cite{bezrukov2018new}}
 The graph $(K_9-C_9)^n$ is defined recursively as the cartesian product of  $(K_9-C_9)$ with $(K_9-C_9)^{(n-1)}$, $n\geq 2$, where  $(K_9-C_9)$ is the graph obtained by deleting a cycle on 9 vertices from the complete graph $K_9$ on 9 vertices.    
\end{definition}
\section{Main Results}
In this section, we prove that $(K_9-C_9 )^n$, $n\geq 1$, is pancyclic. $(K_9-C_9 )^n$ comprises of underlying meshes $M_1,M_2,...M_{9^{n-2}}$, as subgraphs, each isomorphic to the mesh $M_{9 \times 9}$ with 9 rows and 9 columns. Let $(i,j)_k$ denote the vertex at the $i^{th}$ row and $j^{th}$ column of $M_k$, $1 \leq i,j \leq 9$, $1 \leq k \leq 9^{n-2}$. Without loss of generality, we write $(i,j)_1$ as $(i,j)$. By Theorem 3, $(K_9-C_9 )$ is pancyclic as it is on 9 vertices and 27 edges and is not a complete bipartite graph. See Figure 1 (c). \newline To prove the pancyclicity of $(K_9-C_9 )^n$, $n\geq 2$, we give Even Pancyclicity Algorithm to list out cycles of even length between 4 and $9^{n}-1$ in $(K_9-C_9 )^n$ and Odd Pancyclicity Algorithm to list out cycles of odd length between 3 and $9^n$. We begin with $n=2$ and extend it to arbitrarily integer $n$.The graph $(K_9-C_9)^n$, $n\geq 1$  has $9^n$ vertices and $3n\times 9^n$ edges. It is $6n-$regular.\\ \\
{\bf Algorithm: Even Pancyclicity of $(K_9-C_9)^2$}\\ \\
{\bf Input:} Cartesian Product graph $(K_9-C_9 )^2$ \\
{\bf Algorithm:} \\
The vertices $(1,1)$, $(1,2)$, $(2,2)$ and $(2,1)$ induce a 4-cycle.\\
{\bf Step 1}  Successively removing the edge $((i,1),(i,2))$ and adding the vertices $(i+1,1)$ and $(i+1,2)$ induce cycles of length $2(i+1)$ at every stage, $2 \leq i \leq 8$.\\
{\bf Step 2} For $1 \leq i \leq 4$, successively removing the edge $((2i-1,j),(2i,j))$ and adding the vertices $(2i-1,j+1)$ and $(2i,j+1)$, induce cycles of length $14i+2j+2$, $2 \leq j \leq 8$ at every stage.\\
{\bf Step 3} For $1 \leq i \leq 3$, successively removing pairs of edges $((8,2i+1),(8,2i+2))$ and adding the vertices $(9,2i+1)$ and $(9,2i+2)$ generates cycles of length $74+2i$.\\ 
{\bf Output}: Even pancyclicity in $(K_9-C_9 )^2$ with cycles of lengths $4, 6, 8,...,80$.\\
See Figure 2. For convenience $(i,j)$ has been written as $ij$
\begin{figure}[h]
    \centering
    \includegraphics[width=15cm]{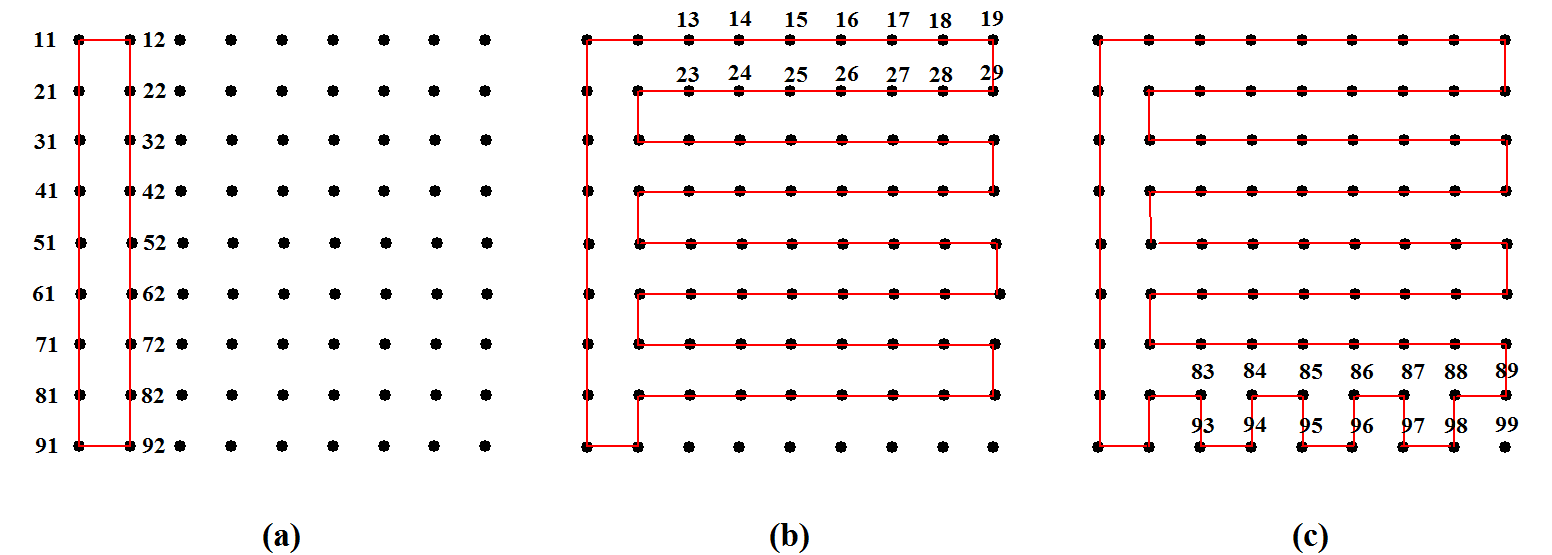}
    \caption{Illustrating even cycles in $M_1$ of length (a) $2(i+1)$, $2 \leq i \leq 8$ (b) $14i+2j+2$, $1 \leq i \leq 4$, $2 \leq j \leq 8$ (c) $74+2i$, $1 \leq i \leq 3$}
    \label{Figure 2}
\end{figure} 
{\bf Proof of correctness :}
 Step 1 begins with a cycle of length 4 and at every stage, the length of already generated even cycle gets incremented by 2 because for $2 \leq i \leq 8$, $2(i+1+1)-2(i+1)=2$. Thus the cycles generated in Step 1 are $6,8,\dots,18$.\\
In Step 2, for a fixed $i$, $1\leq i \leq 4$, $(14i+2(j+1)+2)-(14i+2j+2)=2$, for $2\leq j \leq 8$. Thus for $1 \leq i \leq 4$, the length of cycles generated are $14i+2j+4$, $1 \leq j \leq 7$. Step 3  generates cycles of length $76,78$ and $80$.\\ \\
{\bf Algorithm: Odd Pancyclicity of $(K_9-C_9)^2$}\\ \\
{\bf Input:} Cartesian Product $(K_9-C_9 )^2$ \\
{\bf Algorithm:}\\ 
{\bf Step 1} For $1 \leq j \leq 3$, the vertices $(1,1),(2,1),....,(2j+1,1)$ together with the edge $((1,1),(2j+1,1))$ generates cycles of length $2j+1$.\\ 
{\bf Step 2} (a) For $1\leq i \leq 3$, successively removing the edge $((2i-1,j),(2i,j))$ and adding the vertices $(2i-1,j+1)$ and $(2i,j+1)$ for $1 \leq j \leq 8$, generates cycles of length $16i+2j-9$ at every stage.\\
(b) Deleting edge $((6,1),(7,1))$ and $((6,1),(6,2))$ and adding the vertices $(8,1), (8,2)$ and $(7,2)$ generates a cycle of length 57.\\
(c) Successively removing the edge $((7,i),(8,i))$ and adding the vertices $(7,i+1)$ and $(8,i+1)$, $2\leq i \leq 8$ generates cycles of length $55+2i$ at every stage.\\
(d) For $1 \leq i \leq 4$, successively deleting the edge $((8,2i-1),(8,2i))$ and adding the vertices $(9,2i-1)$ and $(9,2i)$ generates cycles of length $71+2i$.\\
(e) Removing the edges $((4,1),(5,1))$, $((9,3),(9,4))$ and adding edges $((4,1),(6,1))$, $((5,1),(6,1))$, $((9,3),(9,9))$ and $((9,4),(9,9))$ generates a cycle of length 81.\\
{\bf Output:} Odd Pancyclicity in $(K_9-C_9 )^2$ with cycles of length $3, 5, 7,\dots,81$.\\
\begin{figure}[h]
    \centering
    \includegraphics[width=15cm]{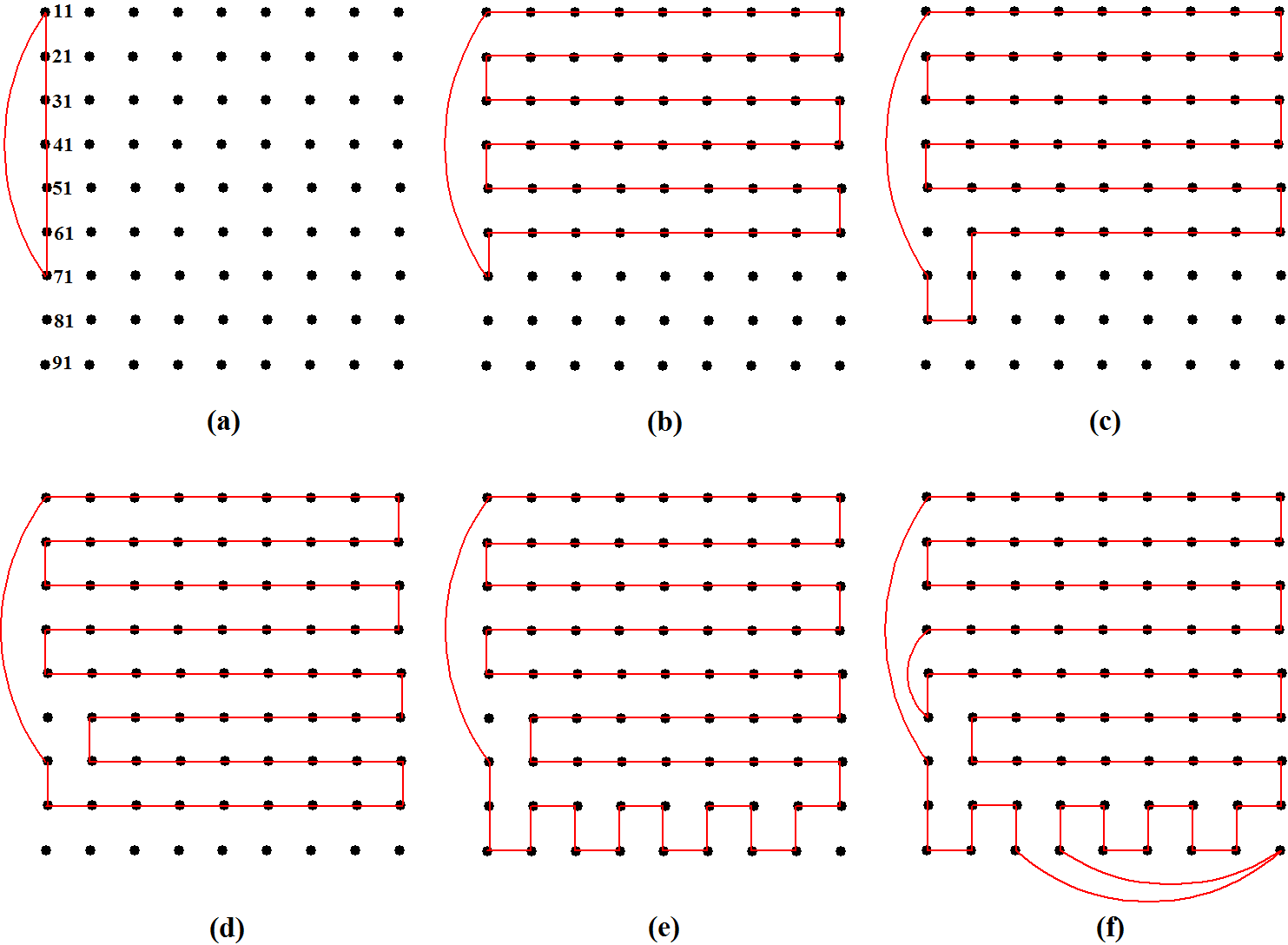}
    \caption{Illustrating odd cycles in $M_1$ of length (a) $2j+1$, $1 \leq j \leq 3$ (b) $16i+2j-9$, $1\leq i \leq 3$, $1 \leq j \leq 8$ (c) 57 (d) $55+2i$, $2\leq i \leq 8$ (e) $71+2i$, $1 \leq j \leq 4$ (f) 81}
    \label{Figure 3}
\end{figure}
{\bf Proof of correctness :}
 Step 1 generates cycles of length $3, 5$ and $7$. Step 2 (a) generates cycles of length $9, 11,\dots,55$.  Step 2 (b) generates cycles of length $57$, Step 2 (c) generates cycles of length $59, 61,\dots,71$ and Step 2 (d) generates cycles of length $73, 75,\dots,79$. Step 2 (e) generates a cycle of length $81$.\\  \\
{\bf Algorithm: Pancyclicity of $(K_9-C_9)^n$}\\ \\
By the Structure of $(K_9-C_9)^n$, the meshes $M_1, M_2,...,M_{9^{n-2}}$ are connected sequentially with edges between pairs of vertices as follows:
(a) $(1,2)_{k-1}$ and $(1,2)_k$\\
(b) $(1,3)_{k-1}$ and $(1,1)_k$, for $2 \leq k \leq 9^{n-2}$.\\ 
{\bf Input:} Cartesian Product graph $(K_9-C_9)^n$ \\
{\bf Algorithm:}\\
{\bf Step 1} Execute Algorithm Even Pancyclicity of $(K_9-C_9)^2$ in $M_1$. Delete the edge $((1,2)_1,(1,3)_1)$ and add the edges $((1,2)_1,(1,2)_2)$, $((1,3)_1,(1,1)_2)$ and $((1,1)_2,(1,2)_2)$ to generate a cycle of length $82$. Execute Algorithm Even Pancyclicity of $(K_9-C_9)^2$ in $M_2$. Inductively remove the edge $((1,2)_{k-1},(1,3)_{k-1})$ and add the edges $((1,2)_{k-1},(1,2)_{k})$, $((1,3)_{k-1},(1,1)_{k})$ and $((1,1)_k,(1,2)_k)$ and execute algorithm even pancyclicity in $M_{k}$, $2 \leq k \leq 9^{n-2}$.\\
{\bf Step 2} Execute Algorithm odd Pancyclicity of $(K_9-C_9)^2$ in $M_1$. Delete the edge  $((1,2)_1,(1,3)_1)$ and successively add pairs of vertices (i) $(1,1)_2, (2,1)_2$ (ii) $(3,1)_2, (4,1)_2$ (iii) $(5,1)_2, (6,1)_2$  generating cycles of length 83, 85 and 87 at every stage. Repeat the procedure for $k$, $2 \leq k \leq 9^{n-2}$ as follows: execute the Algorithm Odd Pancyclicity of $(K_9-C_9)^2$ in $M_k$.  Inductively remove the edge $((1,2)_{k-1},(1,3)_{k-1})$ and add the edges $((1,2)_{k-1},(1,2)_{k})$ and $((1,3)_{k-1},(1,1)_{k})$ and execute Algorithm Odd Pancyclicity in $M_{k}$, $2 \leq k \leq 9^{n-2}$.\\
\begin{figure}[h]
    \centering
    \includegraphics[width=15cm]{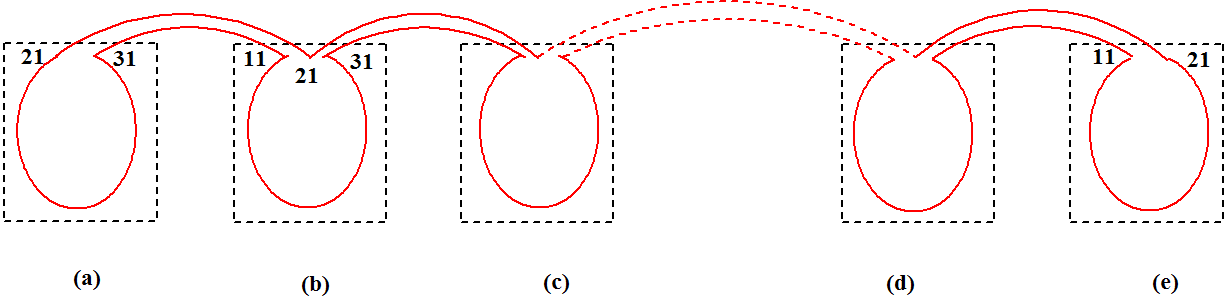}
    \caption{Constructing cycles in (a) $M_1$ (b) $M_2$ (c) $M_3$ (d) $M_{9^{n-2}}-1$ (e) $M_{9^{n-2}}$ for $n \geq 3$}
    \label{Figure 4}
\end{figure}
{\bf Step 3} In $(K_9-C_9)^k$, consider the sequence of meshes $M_1$,$M_2$,\dots,$M_{9^{k-2}}$. (i) Execute Algorithm Even Pancyclicity for meshes $M_1$,$M_2$,\dots,$M_{9^{n-2}-1}$ together with Step 1. (ii) Execute Algorithm Odd Pancyclicity for mesh $M_{9^{k-2}}$, $3\leq k\leq 9^{n-1}$, $n\geq 3$ together with Step 2.\\
{\bf Output:} Pancyclicity in $(K_9-C_9)^n$.\\ \\
{\bf Proof of correctness :}
 Step 1 generates cycles of length $82, 84,....,9^n-1-(9^{n-2}-1)$. Step 2 generates cycles of length  $83, 85,..,9^n$. There are $\frac{9^{n-2}-1}{2}$ cycles of even length not accounted for in Steps 1 and 2. Step 3 generates the left out even cycles of length $80(9^{n-2})+j$, $2\leq j \leq 9^{n-2}-1$, where $j$ is even.
\section{Conclusion}
In this paper we have investigated the pancyclicity of $(K_9-C_9)^n$ network which helps in the development of efficient parallel algorithms.
 \bibliographystyle{unsrt} 
 \bibliography{pancycle}
\end{document}